\documentclass[10pt]{amsart}
\usepackage{amsmath, amssymb, amsthm}

\author{Alexander Berenstein}
\address{Alexander Berenstein
\\ University of Illinois at Urbana-Champaign \\ 1409 West
  Green Street\\Urbana, IL 61801-2975\\USA.}
\email{aberenst@math.uiuc.edu}

\thanks{The author would like to thank C. Ward Henson
for usefull comments}

\makeatletter
\def\newrefformat#1#2{%
  \@namedef{pr@#1}##1{#2}}
\def\prettyref#1{\@prettyref#1:}
\def\@prettyref#1:#2:{%
  \expandafter\ifx\csname pr@#1\endcsname\relax%
    \PackageWarning{prettyref}{Reference format #1\space undefined}%
    \ref{#1:#2}%
  \else%
    \csname pr@#1\endcsname{#1:#2}%
  \fi%
}
\makeatother

\newrefformat{eq}{\textup{(\ref{#1})}}
\newrefformat{cha}{Chapter \ref{#1}}
\newrefformat{sec}{Section \ref{#1}}
\newrefformat{tab}{Table \ref{#1} on page \pageref{#1}}
\newrefformat{fig}{Figure \ref{#1} on page \pageref{#1}}

\makeatletter
\def\indsym#1#2{%
  \setbox0=\hbox{$\m@th#1x$}%
  \kern\wd0%
  \hbox to 0pt{\hss$\m@th#1\mid$\hbox to 0pt{$\m@th#1^{#2}$}\hss}%
  \lower.9\ht0\hbox to 0pt{\hss$\m@th#1\smile$\hss}%
  \kern\wd0} 
\def\nindsym#1#2{%
  \setbox0=\hbox{$\m@th#1x$}%
  \kern\wd0%
  \hbox to 0pt{\mathchardef\nn="3236\hss$\m@th#1\nn$\kern1.4\wd0\hss}
  \hbox to 0pt{\hss$\m@th#1\mid$\hbox to 0pt{$\m@th#1^{#2}$}\hss}%
  \lower.9\ht0\hbox to 0pt{\hss$\m@th#1\smile$\hss}%
  \kern\wd0}

\makeatother

 \def\bx{\bar x} \def\by{\bar y}

  \def\U{\mathcal{U}}

 \def\A{\mathcal{A}}

\def\bra{\langle} \def\ket{\rangle} \def\vphi{\varphi}

\def\H{\mathcal{H}}
\def\U{\mathcal{U}}
\def\Oplus{\oplus}

\newcommand{\mynewthm}[3][dummythm]{%
  \newtheorem{#2}[#1]{#3}%
  \newrefformat{#2}{#3 \ref{##1}}%
}

\swapnumbers

\theoremstyle{plain}
\mynewthm{theo}{Theorem}
\mynewthm{prop}{Proposition}
\mynewthm{lema}{Lemma}
\mynewthm{fact}{Fact}
\mynewthm{coro}{Corollary}
\mynewthm{ques}{Question}
\mynewthm{exer}{Exercise}

\theoremstyle{definition}
\mynewthm{defi}{Definition}
\mynewthm{rema}{Remark}
\mynewthm{obse}{Observation}
\mynewthm{conv}{Convention}
\mynewthm{nota}{Notation}
\mynewthm{exam}{Example}

\DeclareMathOperator{\dcl}{dcl}

\title{Hilbert spaces with generic groups of automorphisms}

\begin{document}
\maketitle
\normalsize

\begin{abstract}
Let $G$ be a countable group. We proof that there is a model companion for the approximate theory of a Hilbert space with a group $G$ of automorphisms. We show that $G$ is amenable if and only if the structure induced by countable copies of the regular representation of $G$ is existentially closed.
\end{abstract}

\section{Introduction}
This paper deals with representations of groups on the collection of unitary maps $\U$ of a Hilbert space $\H$. The aim of this note is to identify the collection of such representations that are rich from the model theoretic perspective, that is, the ones that induce existentially closed expansions of the underlying  Hilbert space. 

We will use approximate semantics \cite{HI} to study the expansions of Hilbert spaces by a group of automorphisms. This approach to model theory falls into the larger setting of compact abstract theories \cite{BY0,BY1}.

There are several papers that deal with expansions of Hilbert spaces with bounded linear operators. Iovino and Henson \cite{Io1} proved that the approximate theory of each of these structures is stable. Their argument shows that there is a bound on the density character of the space of types. 

In \cite{BB} Berenstein and Buechler study Hilbert spaces expanded by a commuting family of normal operators. The expansion of a Hilbert space with an abelian group of automorphisms is an example of this kind of structures. Their approach to stability is via simplicity and the study of non-dividing, a combinatorial notion that induces a dimension theory for elements in the Hilbert space. It is proved in \cite{BB} that the approximate theory of Hilbert spaces expanded by a commuting family of normal operators has quantifier elimination and there is a characterization of non-dividing in terms of the spectral decomposition. The key tool in this work was the Spectral Theorem \cite{AG}.

This paper arose from the following problem. Henson and Berenstein have been studying the expansions of the measure algebra of a probability space by a group of automorphisms. They showed the existence of a model companion when the group involved is amenable. It is still an open question whether there is a model companion for an arbitrary countable group $G$. The authors could show that the regular representation $2^G$ with the left $G$-action and the Haar measure gives rise to an existentially closed structure if and only if $G$ is amenable. Henson asked if there was always a model companion for the approximate theory of probability spaces with a group $G$ of automorphisms.

In this paper, we give a positive answer to Henson's question.

Amenability comes into the picture in a similar way as it did above. Countable copies of the regular representation $l^2(G)$ of a group $G$ with the left $G$-action give rise to an existentially closed structure if and only if $G$ is amenable. The key ingredient used in the proof is a Theorem by Hulanicki (see Proposition \ref{saturation}).

This paper is organized in the following way. In section 2 we deal with the existence of a model companion. We also show that the class of existentially closed structures has the amalgamation property and thus it has quantifier elimination. In section 3 we will give a characterization of non-dividing in these theories and show a weak form of elimination of imaginaries. We also prove these structures are superstable.

\section{Generic automorphisms of Hilbert spaces}

Let $G$ be a discrete group. To simplify the presentation of the 
arguments, we will assume that $G$ is countable, but the arguments
 can be easily generalized to the uncountable case. The main result
 about unitary representations that we will use in this paper 
is a generalizations of a theorem of Hulanicki that can be found
 in \cite{Zi}. We will follow 
the presentation and notation from \cite{HK,Ke} and
concentrate our attention on separable Hilbert spaces, where the
Hilbert space has an infinite countable orthonormal basis. 
This amounts to studying the separable models of the approximate 
theory of Hilbert spaces expanded by a group $G$ of automorphisms.

We start with some basic definitions that we copied from \cite{Ke}.

Let $G$ be a countable group. A \emph{unitary representation} 
of $G$ on a separable Hilbert space $H$ is an action of $G$ on $H$ by 
unitary operators, or equivalently, a homomorphism $\pi: G\to U(H)$, 
where $U(H)$ is the unitary group of $H$. The action of $g\in G$ on
 $v\in H$ according to $\pi$ is denoted by $\pi(g)(v)$ or simply $g(v)$.
We also denote $H$ by $H_\pi$. We denote the representation by $(H_\pi,\pi)$ or by $\pi$.

Given a countable family $\{\pi_i:i\in I\}$ of unitary representations 
of $G$ with $H_{i}=H_{\pi_i}$, we define their \emph{direct sum} $\Oplus_i \pi_i$ as follows: Let $H=\Oplus_i H_i$ and 
$g(\oplus_i v_i)= \oplus \pi_i(g)(v_i)$. We denote, for each unitary 
representation $\pi$, by $n\pi$ the direct sum of $n$ copies of $\pi$, where $1\leq n\leq \aleph_0$. We also write $\infty \pi$ for $\aleph_0\pi$.

\begin{defi}
Let $\pi$, $\rho$ be unitary representations of $G$. Let $\epsilon>0$ 
and $F\subset G$ finite. We say that $\rho$ is 
\emph{$(\epsilon,F)$-contained} in $\pi$ if for every 
$v_1,\dots,v_n\in H_{\rho}$, there are $w_1,\dots,w_n\in H_\pi$ such 
that $|\bra \pi(g)v_i,v_j \ket-\bra  \rho(g)w_i),w_j\ket |<\epsilon$ for 
all $g\in F$. We say that $\rho$ is 
\emph{weakly contained in $\pi$} and write $\rho \prec \pi$ if $\rho$ 
is $(\epsilon,F)$ contained in $\pi$ for every $\epsilon$, $F$.
\end{defi}

For each countable group $G$, we denote by $\lambda_G$ the 
\emph{regular representation} of $G$, that is, the left-shift action 
of $G$ on $l^2(G)$ defined by $\lambda_G(g) f(h)=f(g^{-1}h)$. We 
denote by $1_G$ the trivial representation, the one that sends every 
element in $G$ to the identity map. 

We denote by $L$ the language of Hilbert spaces and by $L_G$ the language $L$ expanded by unary functions $\pi(g)$ for each $g\in G$

When we study the model obtained by the expansion of a Hilbert space 
with a representation $\pi$, we write the structure as 
$(H_\pi,+,\bra,\ket,\pi(g):g\in G)$. We add a second sort for the 
complex field and constants for all rational numbers. In the language 
for these structures, atomic formulas are of the form $\bra s,t \ket \leq p$ 
\ $\bra s,t \ket 
\geq q$, where $s,t$ are terms in the sort of the Hilbert space
and $p,q$ are non-negative rational numbers. 
The collection of positive formulas is built by 
closing the set formed by the atomic formulas under conjunctions, 
disjunctions and bounded existential and universal quantifiers. \emph{Approximations} 
to a given formula are constructed inductively by relaxing the 
conditions on the atomic formulas. For example, an approximation 
of the formula $\bra s,t \ket\leq p$ is a formula of the 
form $\bra s,t\ket \leq p+\epsilon$, where $\epsilon$ is a positive 
rational number. A more detailed explanation of positive formulas and
their approximations can be found in \cite{HI}.

We denote by $T$ the approximate theory of Hilbert spaces and by $T_G$ the union of $T$ with the axioms stating that each $\tau_g$ is an automorphisms and that $\tau_{g}:g\in G$ is a representation of $G$. We denote by $tp$ the types in the language $L$ and by $tp_G$ the types in the language $L_G$. When $G$ is clear from context and $(H_\pi,\pi)$ is a representation of $G$, we write $(H_\pi,+,\bra,\ket,\pi)$ instead of $(H_\pi,+,\bra,\ket,\pi(g):g\in G)$

\begin{defi}
Let $(H_\pi,\pi)$ be a representation. We say that 
$(H_\pi,+,\bra,\ket,\pi)$ is \emph{existentially
closed} if given any representation $(H_\eta,\eta)\supset (H_\pi,\pi)$, 
elements $v_1,\dots,v_n\in H_\pi$, 
quantifier free formula $\vphi(\bx,\by)$ such that 
$(H_\eta,+,\bra,\ket,\eta)\models \exists \bx \vphi(\bx,v_1,\dots ,v_n)$ 
and an approximation $\vphi'(\bx,\by)$ of $\vphi(\bx,\by)$,
$(H_\pi,+,\bra,\ket,\pi)\models \exists \bx \vphi'(\bx,v_1,\dots ,v_n)$.
\end{defi}

The first goal of this section
 is to show that $(H_{\infty  \lambda_G},+,\bra,\ket,\infty\lambda_G)$ 
is existentially closed if and only if $G$ is amenable. We will need the
 following generalization of Hulanicki's Theorem:

\begin{prop}(Proposition 7.3.6 in \cite{Zi}))\label{saturation}
(a) Suppose that $G$ is amenable. Then for any representation $\rho$ of $G$, $\rho\prec \infty \lambda_G$.\\
(b) If $1_G\prec \infty\lambda_G$, then $G$ is amenable.
\end{prop}

\begin{coro} Assume that $G$ is not amenable.
Then the structure $(H_{\infty\lambda_G},+,\bra,\ket,\infty \lambda_G)$ 
is not existentially closed.
\end{coro}

\begin{proof}
let $G$ be non-amenable, let $(H_{1_G},1_G)$ be the trivial representation and assume, as a way to a contradiction, that $(H_{\infty\lambda_G},+,\bra,\ket,\infty \lambda_G)$ is existentially closed. Then $(H_{\infty \lambda_G}\oplus H_{1_G},+,\bra,\ket,\infty \lambda_G \oplus 1_G)$ is an extension of $(H_{\infty\lambda_G},+,\bra,\ket,\infty \lambda_G)$ and thus $1_G\prec \infty \lambda_G$. This contradicts Proposition \ref{saturation}.
\end{proof}

\begin{theo}\label{superstructure} Assume that $G$ is amenable 
and let $(H_\pi,\pi)$ be any representation. Then the structure 
$(H_\pi \oplus H_{\infty\lambda_G},+,\bra,\ket,\pi \oplus \infty \lambda_G)$
is existentially closed.
\end{theo}

\begin{proof}
Write $(H_\eta,+,\bra,\ket,\eta)$ for $(H_\pi \oplus H_{\infty\lambda_G},+,\bra,\ket,\pi \oplus \infty \lambda_G)$. To prove the Theorem it is enough 
to show that any existential formula with parameters in $H_\eta$
 realizable in a separable
extension of $(H_\eta,+,\bra,\ket,\eta)$ is already approximately 
realized in the structure. So let $(H_\rho,+,\bra,\ket,\rho)$ be a 
separable superstructure of $(H_\eta,+,\bra,\ket,\eta)$.
Let $v_1,\dots,v_n\in H_{\eta}$ and let $v_{n+1},\dots,v_{n+m}\in 
H_\rho$. Let $v_i'=v_i$ for $i\leq n$.

\textbf{Claim} For any $\epsilon>0$ and $F\subset G$ finite, there 
are $v_{n+1}',\dots,v_{n+m}'\in 
H_{\eta}$ such that $|\bra \rho(g)v_i,v_j \ket-
\bra \eta(g)v_i',v_j'\ket |<\epsilon$ for $i,j\leq n+m$.

By the perturbation Lemma in \cite{HI}, it is enough to show the claim 
when the parameters $v_1,\dots,v_n$ come from a dense subset of $H_\pi\oplus H_{\infty \lambda_G}$. Identify
$H_\pi \oplus H_{i\lambda_G}$ with its canonical embedding in 
$H_\pi \oplus H_{\infty \lambda_G}$. Then 
$\cup_{i\in \omega}H_\pi\oplus H_{i \lambda_G}$ is dense in 
$H_\pi \oplus H_{\infty \lambda_G}$ and we may assume that 
$v_1,\dots,v_n\in H_\pi \oplus H_{l \lambda_G}$ for some $l$. 
Since $H_\pi\oplus H_{l \lambda_G}\subset H_\pi \oplus H_{\infty \lambda_G} \subset H_\rho$, we can write $v_i=u_{i}+w_{i}$, where $u_{i}$ in the projection of $v_i$ 
on $H_\pi \oplus H_{l\lambda_G}$. Note that $v_i=u_i$ for $i\leq n$.

Let $w_i'=0$ for $i\leq n$. By Proposition \ref{saturation} 
there are $w_{n+1}',\dots,w_{n+m}'\in 
H_{\infty \lambda_G}\cap H_{l \lambda_G}^\perp \cong H_{\infty \lambda_G}$ such that $|\bra \rho(g) w_i,w_j \ket-
\bra \infty \lambda_G(g)w_i',w_j'\ket |<\epsilon$ for 
$n+1\leq i,j\leq n+m$. Define $v_i'=u_i+w_i'$ for $i\leq n+m$. 
Then $|\bra \rho(g) v_i,v_j \ket- \bra \infty 
\lambda_G(g)v_i',v_j'\ket |<\epsilon$ for $i,j\leq n+m$.
\end{proof}

Let $(H_\pi,\pi)$, be a representation and let $(H_\rho,\rho)$, 
$(H_\eta,\eta)$ be representations extending $(H_\pi,\pi)$. That is,
$(H_\pi,+,\bra,\ket ,\pi)\subset (H_\rho,+,\bra,\ket ,\pi)$ and
$(H_\pi,+,\bra,\ket ,\pi)\subset (H_\rho,+,\bra,\ket ,\pi)$. Then 
there are representations $(H_{\rho'},\rho')$ and $(H_{\eta'},\eta')$ 
such that $(H_\rho,\rho)=(H_\pi \oplus H_{\rho'},\pi\oplus \rho')$ and 
$(H_\eta,\eta)=(H_\pi\oplus H_{\eta'},\pi \oplus \eta')$. The 
\emph{amalgamation} of $(H_\rho,\rho)$  and $(H_\eta,\eta)$ over 
$(H_\pi,\pi)$ is $(H_\pi\oplus H_{\eta'}\oplus H_{\rho'},\pi \oplus \eta'\oplus \rho')$. This proves that the class
of Hilbert spaces expanded by a group $G$ of automorphisms has the
\emph{amalgamation property}.

Let $\Sigma_G$ be the 
collection of existential formulas that are true for
$(H_{\infty \lambda_G},+,\bra,\ket, \infty \lambda_G)$. Let $\Sigma^{-}_G$ 
be the collection of approximations of formulas in $\Sigma_G$. We will prove 
that when $G$ is amenable, $\Sigma^{-}_G\cup T_G$ is an axiomatization for the
 class of existentially closed expansions of a Hilbert space with a group
 $G$ of automorphisms. 

\begin{prop}\label{p1}
Let $G$ be amenable and let $(H_\pi,\pi)$ be a representation of $G$ such that $(H_\pi,+,\bra,\ket,\pi)\models \Sigma_G^-$. 
Then  $(H_\pi,+,\bra,\ket,\pi)$ is existentially closed.
\end{prop}

\begin{proof}
Let $(H_\pi,+,\bra,\ket,\pi)\models \Sigma^-_G$. Then there is a separable
 approximate elementary superstructure $(H_\rho,+,\bra,\ket,\rho)$ of 
$(H_\pi,+,\bra,\ket,\pi)$ that contains $(H_{\infty \lambda_G},+,\bra,\ket, \infty \lambda_G)$ as a 
substructure. By Theorem \ref{superstructure}, $(H_\rho,+,\bra,\ket,\rho)$
 is existentially closed. 

Now let $(H_\eta,+,\bra,\ket,\eta) \supset (H_\pi,+,\bra,\ket,\pi)$. Then we can amalgamate $(H_\eta,+,\bra,\ket,\eta)$ and $(H_\rho,+,\bra,\ket,\rho)$
over $(H_\pi,+,\bra,\ket,\pi)$. Since $(H_\rho,+,\bra,\ket,\rho)$ is
 existentially closed, any existential formula with parameters in $H_\pi$ 
true in $(H_\eta,+,\bra,\ket,\eta)$ is approximately
true in $(H_\rho,+,\bra,\ket,\rho)$ and thus it is also approximately true in 
its approximately elementary substructure $(H_\pi,+,\bra,\ket,\pi)$.
\end{proof}

\begin{prop}\label{p2} Let $G$ be amenable and let $(H,+,\bra,\ket,\pi)$ be 
existentially closed. Then $(H_\pi,+,\bra,\ket,\pi)\models \Sigma^-_G$.
\end{prop}

\begin{proof}
Let $(H_\pi,+,\bra,\ket,\pi)$ be existentially closed. Then
$(H_{\infty\lambda_G}\oplus H_\pi,+,\bra,\ket,\infty \lambda_G \oplus \pi)$
 is a superstructure of $(H_\pi,+,\bra,\ket,\pi)$ that models $\Sigma_G^-$ and thus $(H_\pi,+,\bra,\ket,\pi)\models \Sigma^-_G$.
\end{proof}

\begin{theo} Let $G$ be amenable and let $T_{GA}=Th_{\A}(H_{\infty\lambda_G},+,\bra,\ket,\infty \lambda_G)$. Then $T_{GA}$ has quantifier elimination and it is 
axiomatized by $T_G \cup \Sigma_G^-$.
\end{theo}

\begin{proof}
The models of $\Sigma_G^-$ are existentially closed and can 
be amalgamated, thus all of them have the same approximate theory. 
Furthermore, using the amalgamation property and the fact that
$T_{GA}$ is model complete, it is easy to show that the theory 
has quantifier elimination. Finally the claim about axiomatizability 
follows from Propositions \ref{p1} and \ref{p2}.
\end{proof}

\begin{coro} Let $G$ be an amenable group. Let $(H_\pi,\pi)$ be 
a representation. Then $(H_\pi,+,\bra,\ket,\pi)$ 
existentially closed if and only if for all representations 
$(H_\sigma,\sigma)$ we have $\sigma\prec \pi$.
\end{coro}

Ben-Yaacov proved that when $G=\mathbb{Z}$, $(H_\pi,+,0,\bra,\ket,\pi)$ is existentially closed if and only if the spectrum of $\pi(1)$ is $S^1$. This result was also proved independently by Usvyatsov and Zadka. It is straightforward to generalize their result result to abelian groups.

\begin{ques}
Let $G$ be amenable. Is there a spectral characterization of 
all representations $\pi$ such that $\pi\prec \infty\lambda_G$ and 
$\infty\lambda_G \prec \pi$? 
\end{ques}

\begin{rema}
The condition of amenability of $G$ was not essential to show that
there is an axiomatization for the existentially closed expansions 
of Hilbert spaces with automorphisms. It was sufficient to start 
with an existentially closed model (when $G$ is amenable 
$(H_{\infty\lambda_G},+,\bra,\ket,\infty \lambda_G)$ is existentially
 closed) and then show that its approximate existential theory axiomatizes 
the class of existentially closed models. Finally observe that for any 
countable group $G$ there is a representation $\pi$ such that $(H_\pi,+,\bra,
\ket,\pi)$ is existentially closed.
\end{rema}

\newpage

\section{Stability}

Let $G$ be a countable group, let $(H_\eta,+,\bra,\ket,\eta)$ be existentially closed and let $T_{GA}=Th_{\A}(H_\eta,+,\bra,\ket,\eta)$. Let $\kappa$ be a cardinal larger than $2^{\aleph_0}$ and let $(H_\pi,+,\bra,\ket,\pi)\models T_{GA}$ be $\kappa$-saturated. By the work of Henson and Iovino \cite{Io1}, the theory $T_{GA}$ is stable. In this section we will prove that $T_{GA}$ is superstable and we will chraracterize non-dividing. We will assume that the reader is familiar with the notions of \emph{definable closure} and \emph{non-dividing}. The reader can check \cite{BB,BH,BY0,BY1} for the definitions.

Let $\dcl$ stand for the definable closure in the language $L$ and let $\dcl_G$ stand for the definable closure in the language $L_G$. It is easy to see that for $A\subset H_\pi$ such that $|A|<\kappa$, $\dcl_G(A)=\dcl(\cup_{g\in G, a\in A}\pi(g)(a))$.

Let $B,C\subset H_\pi$ be of cardinality less than $\kappa$,
let $(a_1,\dots,a_n) \in H^n$ and assume that $C=\dcl(C)$, 
so $C$ is a Hilbert subspace of $H_\pi$. Denote by $P_C$ the projection 
on $C$. It is proved
in \cite{BB} that $tp(a_1,\dots,a_n/C \cup B)$ does not divide over $C$ 
if and only if for all $i\leq n$ and all $b\in B$, $a_i-P_C(a_i) \perp b-P_C(b)$. In particular, non-dividing is \emph{trivial}. 
The set $\{P_C(a_1):i\leq n\}$ forms the smallest set of elements 
in $C$ over which $tp(a_1,\dots,a_n/C)$ does not divide. Such a set is called a \emph{built-in canonical base} for $tp(a_1,\dots,a_n/C)$ (see \cite{BB,BH} for a formal definition).

Let $(g_j:j\in \omega)$ be an enumeration of $G$. A standard argument shows that $tp_G(a_1,\dots,a_n/B\cup C)$ does not $L_G$-divide over $C$ if and only if $tp(\dcl_G(a_1,\dots,a_n)/\dcl_G(B\cup C))$ does not divide over $\dcl_G(C)$ if and only if $tp(\pi(g_j)(a_i)/\dcl_G(B\cup C))$ does not divide over $\dcl_G(C)$ for all $j\in \omega$, $i\leq n$. This gives a characterization 
of non-dividing for $T_{GA}$. Note that non-diving is also trivial in $(H_\pi,+,0,\bra,\ket,\pi)$. From this characterization of non-dividing we can conclude that $T_{GA}$ is stable and that types over sets are stationary.

Assume now that $C=\dcl_G(C)$ and let $a_1,\dots,a_n\in H_\pi$. Then $tp_G(a_1,\dots,a_n/C)$ does not divide over $\{P_C(a_1),\dots,P_C(a_n)\}$. In particular, up to interdefinability, $\{P_C(a_1),\dots,P_C(a_n)\}$ is the smallest subset of $C$ over which $tp_G(a_1,\dots,a_n/C)$ does not divide. This shows the structure $(H_\pi,+,\bra,\ket,\pi)$ has built-in canonical bases.

\begin{defi}
Let $G$ be a countable group and let $(H_\pi,+,\bra,\ket,\pi)\models T_{GA}$ be $\kappa$-saturated.. We say that $(H_\pi,+,\bra,\ket,\pi)$ is \emph{superstable} if for any $a_1,\dots,a_n\in H_\pi$, $A\subset H_{\pi}$ and $\epsilon>0$, there are $b_1,\dots,b_n\in H_\pi$ and $A_0\subset A$ finite such that $\|a_i-b_i\|<\epsilon$ and $tp(b_1,\dots,b_n/A)$ does not divide over $A_0$.
\end{defi}

It was known to Ben-Yaacov, Usvyatsov and Zadka that $T_{GA}$ is superstable when $G=\mathbb{Z}$. We now prove that $T_{GA}$ is superstable for all $G$:

\begin{prop}
Let $G$ be a countable group and let $(H_\pi,+,\bra,\ket,\pi)\models T_{GA}$. Then $(H_\pi,+,\bra,\ket,\pi)$ is superstable.
\end{prop}

\begin{proof}
Let $a_1,\dots,a_n\in H_\pi$, $\epsilon>0$ and $A\subset H_{\pi}$ be of cardinality less than $\kappa$. Let $C=\dcl_G(A)$. Then for any $\epsilon>0$ there is $A_0 \subset A$ finite such that $\|P_C(a_i)-P_{C_{0}}(a_i)\|<\epsilon$ for $i\leq n$, where $C_0=\dcl(A_0)$. Let $b_i=a_i-P_{C}(a_i)+P_{C_0}(a_i)$. Then $tp(b_1,\dots,b_n/A)$ does not divide over $A_0$ and $\|a_i-b_i\|<\epsilon$ for $i\leq n$.
\end{proof}

\end {document}